\documentclass[preprint,11pt]{elsarticle}

\usepackage{amssymb}
\usepackage{amsmath}
\usepackage{amsfonts}
\usepackage{amssymb}
\usepackage{indentfirst,latexsym,bm}
\usepackage{amsthm}
\usepackage[all]{xy}
\usepackage{color,xcolor}
\newtheorem{theorem}{Theorem}[section]
\newtheorem{lemma}[theorem]{Lemma}
\newtheorem{corollary}[theorem]{Corollary}

\theoremstyle{definition}

\theoremstyle{definition}

\theoremstyle{remark}
\newtheorem{remark}{Remark}[section]

\numberwithin{equation}{section}

\journal{XXX}

\begin{document}

\begin{frontmatter}



\title{The $(n+1)$-centered operator on a Hilbert $C^*$-module}

\author[rvt]{Na Liu}
\ead{liunana0616@163.com}
\author[rvt1]{Qingxiang Xu}
\ead{qingxiang$\_$xu@126.com}
\author[rvt1]{Xiaofeng Zhang}
\ead{xfzhang8103@163.com}
\address[rvt]{College of Mathematics and Information Science, Zhengzhou University of Light Industry, Zhengzhou 450002, PR China}
\address[rvt1]{Department of Mathematics, Shanghai Normal University, Shanghai 200234, PR China}

\begin{abstract}
Let $T$ be an adjointable operator on a Hilbert $C^*$-module such that $T$ has the polar decomposition $T=UT|$. For each natural number $n$, $T$ is called an $(n+1)$-centered operator if $T^k=U^k|T^k|$ is the polar decomposition for $1\le k\le n+1$. This paper initiates  the study of the $(n+1)$-centered operator via the generalized Aluthge transform and the generalized iterative Aluthge transform. Some
new characterizations of the $(n+1)$-centered operator are provided.
\end{abstract}


\end{frontmatter}



\section*{introduction}

Motivated by the weighted shifts, the centered operator is introduced in \cite{Morrel} for a bounded linear operator on a Hilbert space.
Let  $T=U|T|$ be the polar decomposition.  It is proved in  \cite[Theorem~I]{Morrel} that whenever  $T$ is a centered operator, $T^k=U^k|T^k|$
is the polar decomposition for all $k\in\mathbb{N}$, where $\mathbb{N}$ denotes the set of positive integers. The reverse of such a property is also true; see \cite[Theorem~3.2]{Ito-Yamazaki-Yanagida} for the details. To get a deeper understanding of the centered operator, conditions are investigated in \cite{Ito,Liu-Luo-Xu-1} for each  $n\in\mathbb{N}$ such that
\begin{equation}\label{defn of n centered operator}T^k=U^k|T^k|\quad\mbox{is the polar decomposition for $1\le k\le n$}.\end{equation}
An operator $T$ satisfying \eqref{defn of n centered operator} is called an $n$-centered operator \cite[Definition~4.1]{Liu-Luo-Xu-1}.
Hence, an operator is centered if and only if it is $n$-centered for all $n\in\mathbb{N}$. It is remarkable that for each $n\in\mathbb{N}$, an operator $T$ acting on certain Hilbert space can be constructed such that $T$ is $(n+1)$-centered, whereas it fails to be $(n+2)$-centered \cite[Theorem~6.2]{Liu-Luo-Xu-1}. This shows the non-triviality in the study of the $(n+1)$-centered operator.
The $2$-centered operator, usually known as the binormal operator \cite[Section~1]{Campbell-1},
plays a crucial  role in the study of the centered operator  via the Aluthge transform and the iterative Aluthge transform \cite[Section~3]{Ito}. In terms of the Aluthge transform, a characterization of the binormal operator can be found in \cite[Theorem~3.1]{Ito-Yamazaki-Yanagida}, which has been improved in \cite[Theorem~4.14]{Liu-Luo-Xu-1} via the generalized Aluthge transform defined by \eqref{equ:defn of generalized Aluthge transforms}\footnote{\,In this paper, the generalized Aluthge transform is defined for every positive numbers $\alpha$ and $\beta$. Some restrictions are however usually imposed on $\alpha$ and $\beta$  in the literature (see e.g., \cite{Furuta,Moslehian}).}. Some interesting results on the binormality of the iterative Aluthge transform can be found in \cite{Ito,Ito-Yamazaki-Yanagida}.

The purpose of this paper is to give some  new characterizations of the $(n+1)$-centered operator acting on a general Hilbert $C^*$-module. Let $T$ be an adjoint operator such that $T$ has the polar decomposition $T=U|T|$. For each $n\in\mathbb{N}$,
a necessary and sufficient condition is given in \cite[Theorem~4.3]{Liu-Luo-Xu-1} under which $T$ is $(n+1)$-centered (see Lemma~\ref{lem:characterization of (n+1)-centered}). Such a result has been improved in Theorem~\ref{thm:further describe of (n+2)-centered}, where
two positive parameters are employed. Based on \cite[Theorem~4.14]{Liu-Luo-Xu-1}, this paper initiates the study of the $(n+1)$-centered operator via the generalized Aluthge transform and the generalized iterative generalized Aluthge transform. As a result, some non-trivial generalizations of \cite{Ito,Ito-Yamazaki-Yanagida} have been obtained. For instance, it can be deduced from \cite[Theorem~3.6]{Ito} and \cite[Theorem~4.3]{Liu-Luo-Xu-1} that for each $n\ge 2$,  $T$ is $(n+1)$-centered if and only if the iterative Aluthge transform $\widetilde{T}^{(k)}$ is binormal for $k=0,1,\cdots, n-1$. It will be shown in Theorem~\ref{thm:binormality of T-alpha-beta-n}  that the latter condition can be replaced with the weaker one: there exist some positive numbers $\alpha$ and $\beta$ such that the  generalized  iterative Aluthge transform  $\widetilde{T}_{\alpha,\beta}^{(k)}$ is binormal for $0\le k\le n-1$, which in turn implies the binormality of $\widetilde{T}_{\alpha^\prime,\beta^\prime}^{(k)}$  for every $k\in\{0,1,\cdots,n-1\}$ and every positive numbers $\alpha^\prime$ and $\beta^\prime$.

The paper is organized as follows.  Some basic knowledge about Hilbert $C^*$-modules and adjointable operators are provided in Section~\ref{sec:preliminaries}. A new characterization of the $(n+1)$-centered operator is given  in Section~\ref{sec:a new characterization} via the introduction of two positive parameters. In Section~\ref{sec:generalized Aluthge transform}, we focus on the study of the $(n+1)$-centered operator in terms of the generalized Aluthge transform.  In Section~\ref{sec:iterated Aluthge transforms}, we  turn to study the $(n+1)$-centered operator via the generalized iterative Aluthge transform. As an application, a new characterization of the centered operator is provided at the end of this section.

\section{Preliminaries}\label{sec:preliminaries}
Hilbert $C^*$-modules are generalizations of Hilbert spaces by allowing inner products to take values in some $C^{*}$-algebras instead of the complex  field \cite{Lance}. Given  Hilbert $\mathfrak{A}$-modules $H$ and $K$ over a $C^*$-algebra $\mathfrak{A}$, let $\mathcal{L}(H,K)$
be the set of operators $T:H\to K$ for which there is an operator $T^*:K\to
H$ such that $$\langle Tx,y\rangle=\langle x,T^*y\rangle \quad \mbox{for every $x\in H$ and $y\in K$}.$$  We call ${\mathcal
L}(H,K)$ the set of adjointable operators from $H$ to $K$. For every
$T\in \mathcal{L}(H,K)$, its  range and  null space are denoted by
${\mathcal R}(T)$ and ${\mathcal N}(T)$, respectively. In case
$H=K$, $\mathcal{L}(H,H)$ which is abbreviated to $\mathcal{L}(H)$, is a
$C^*$-algebra. Let $\mathcal{L}(H)_+$ denote the set of positive elements
in $\mathcal{L}(H)$. By a projection, we  mean an idempotent and self-adjoint element of certain $C^*$-algebra. Throughout the rest of this paper, $\mathbb{N}$ is the set of positive integers, $\mathfrak{A}$ is a $C^*$-algebra, $E$, $H$ and $K$ are Hilbert $\mathfrak{A}$-modules.

Recall that a closed submodule $M$ of  $H$ is said to be
orthogonally complemented  if $H=M\dotplus M^\bot$, where
$$M^\bot=\big\{x\in H: \langle x,y\rangle=0\quad \mbox{for every}\ y\in
M\big\}.$$
In this case, the projection from $H$ onto $M$ is denoted by $P_M$.

An element $U$ of $\mathcal{L}(H,K)$ is said to be a partial isometry if $U^*U$ is a projection on $H$, or equivalently, $UU^*U=U$  \cite[Lemma~2.1]{XF}.
When  $H=K$, let $[A,B]=AB-BA$ be the commutator of $A,B\in\mathcal{L}(H)$.

For every $T\in\mathcal{L}(H,K)$, let $|T|$ denote the square root of $T^*T$. That is, $|T|=(T^*T)^{\frac12}$ and $|T^*|=(TT^*)^{\frac12}$.
 The polar decomposition of $T$ can be represented as
\begin{equation}\label{equ:two conditions of polar decomposition}T=U|T|\quad \mbox{and}\quad U^*U=P_{\overline{\mathcal{R}(T^*)}},\end{equation}
where $U\in\mathcal{L}(H,K)$ is a partial isometry. From \cite[Lemma~3.6 and Theorem~3.8]{Liu-Luo-Xu}, it is known that
$T$ has the polar decomposition represented by \eqref{equ:two conditions of polar decomposition} if and only if
$\overline{\mathcal{R}(T^*)}$ and $\overline{\mathcal{R}(T)}$ are orthogonally complemented in $H$ and $K$, respectively.
In such case, the polar decomposition of $T^*$ exists and can be represented by
\begin{align}\label{equ:the polar decomposition of T star-pre stage}&T^*=U^*|T^*| \quad \mbox{and}\quad  UU^*=P_{\overline{\mathcal{R}(T)}}.
\end{align}
It is known that for the given $T\in\mathcal{L}(H,K)$, there exists at most a partial isometry $U\in\mathcal{L}(H,K)$  satisfying \eqref{equ:two conditions of polar decomposition} \cite[Lemma~3.9]{Liu-Luo-Xu}, and such a partial isometry $U$ may however fail to be existent \cite[Example~3.15]{Liu-Luo-Xu}. In the case that $T$ has the polar decomposition represented by \eqref{equ:two conditions of polar decomposition},  in what follows we simply say that $T$ has the polar decomposition $T=U|T|$.

To get the main results of this paper, we need the following lemmas, which will be used in the sequel.

\begin{lemma}\label{lem:rang characterization-1} {\rm \cite[Proposition 2.7]{Liu-Luo-Xu}} Let $A\in\mathcal{L}(H,K)$ and $B,C\in\mathcal{L}(E,H)$ be such that $\overline{\mathcal{R}(B)}=\overline{\mathcal{R}(C)}$. Then $\overline{\mathcal{R}(AB)}=\overline{\mathcal{R}(AC)}$.
\end{lemma}

\begin{lemma}\label{lem:Range Closure of Ta and T} {\rm (\cite[Proposition 3.7]{Lance} and \cite[Proposition 2.9]{Liu-Luo-Xu})} Let $T\in\mathcal{L}(H,K)$. Then the following statements are valid:
\begin{enumerate}
\item[{\rm(i)}] $\overline {\mathcal{R}(T^*T)}=\overline{ \mathcal{R}(T^*)}$ and $\overline {\mathcal{R}(TT^*)}=\overline{ \mathcal{R}(T)}$;
\item[{\rm(ii)}] If $T\in \mathcal{L}(H)_+$, then $\overline{\mathcal{R}(T^{\alpha})}=\overline{\mathcal{R}(T)}$ for every $\alpha>0$.
\end{enumerate}
\end{lemma}

\begin{lemma}\label{lem:commutative property extended-3}{\rm \cite[Propositions~2.4 and 2.6]{Liu-Luo-Xu}} Let $S\in \mathcal{L}(H)$ and $T\in\mathcal{L}(H)_+$ be such that $[S,T]=0$. Then the following statements are valid:
\begin{enumerate}
\item[{\rm (i)}] $[S,T^\alpha]=0$ for every $\alpha>0$;
\item[{\rm (ii)}] If in addition $\overline{\mathcal{R}(T)}$ is orthogonally complemented, then $\big[S,  P_{\overline{\mathcal{R}(T)}}\big]=0$.
\end{enumerate}
\end{lemma}

\section{A new  characterization of $(n+1)$-centered operators}\label{sec:a new characterization}
Let $T\in\mathcal{L}(H)$ have the polar decomposition $T=U|T|$, and  let $n\in\mathbb{N}$. Recall that $T$ is said to be $n$-centered if \eqref{defn of n centered operator} is satisfied.  A characterization of the $(n+1)$-centered operator is given in \cite{Liu-Luo-Xu-1} as follows.

\begin{lemma}\label{lem:characterization of (n+1)-centered}{\rm \cite[Theorem 4.3]{Liu-Luo-Xu-1}}  Let $T\in\mathcal{L}(H)$ have the polar decomposition $T=U|T|$. Then for each $n\in\mathbb{N}$, the following statements are equivalent:
\begin{enumerate}
\item[{\rm (i)}] $T$ is $(n+1)$-centered;
\item[{\rm (ii)}] $\big[U^k|T|(U^k)^*, |T|\big]=0$ for all $k=1,2,\cdots, n$.
\end{enumerate}
\end{lemma}
The purpose of this section is to make a generalization of the preceding lemma by introducing two positive  parameters $\alpha$ and $\beta$; see Theorem~\ref{thm:further describe of (n+2)-centered} for the details.

\begin{lemma}\label{lem:property of n-centered-stronger}  Let $T\in\mathcal{L}(H)$ have the polar decomposition $T=U|T|$, and let $n\in\mathbb{N}$ be such that $T$ is $(n+1)$-centered. Then
\begin{equation}\label{equ:property of n-centered-stronger} \big[P_{k}(T^*),|T|\big]=\big[P_{k}(T),|T^*|\big]=0\quad \mbox{for $1\leq k\leq n+1$},
\end{equation}
where
\begin{equation}\label{equ:defn of P n T T-sta}
P_k(T)=U^k(U^k)^*\quad \mbox{and}\quad P_k(T^*)=(U^k)^*U^k.
\end{equation}
\end{lemma}

\begin{proof} It  follows easily from Lemma~\ref{lem:Range Closure of Ta and T} that
 $$\overline{\mathcal{R}(|T|)}=\overline{\mathcal{R}(T^*)}\quad\mbox{and}\quad \overline{\mathcal{R}(|T^*|)}=\overline{\mathcal{R}(T)}.$$ Note that $P_1(T^*)=U^*U$ and $P_1(T)=UU^*$,  so by \eqref{equ:two conditions of polar decomposition} and \eqref{equ:the polar decomposition of T star-pre stage}  we have
\begin{equation*} \big[P_{1}(T^*),|T|\big]=\big[P_{1}(T),|T^*|\big]=0.
\end{equation*}

Given every $k$ with $1\leq k\leq n$, since $T^k=U^k|T^k|$ and $T^{k+1}=U^{k+1}|T^{k+1}|$ are polar decompositions,
by \eqref{equ:two conditions of polar decomposition} both $(U^k)^*$ and $(U^{k+1})^*$ are partial isometries such that $$\overline{\mathcal{R}\big[(T^k)^*\big]}=\mathcal{R}\big[(U^k)^*\big]\quad \mbox{and}\quad
\overline{\mathcal{R}\big[(T^{k+1})^*\big]}=\mathcal{R}\big[(U^{k+1})^*\big],$$
which lead to
$$P_{k+1}(T^*)|T|U^*\cdot(T^k)^*=P_{k+1}(T^*)(T^{k+1})^*=(T^{k+1})^*=|T|U^*\cdot(T^k)^*.$$
Therefore,
$$P_{k+1}(T^*)|T|U^*\xi=|T|U^*\xi,\quad\forall\,\xi\in\overline{\mathcal{R}\big[(T^k)^*\big])}=\mathcal{R}\big[(U^k)^*\big].$$
It follows that
$P_{k+1}(T^*)|T|U^*\cdot (U^k)^*=|T|U^*\cdot (U^k)^*$, which gives
$$P_{k+1}(T^*)|T|P_{k+1}(T^*)=|T|P_{k+1}(T^*)$$ by multiplying $U^{k+1}$ from the right side. Taking $*$-operation on the latter equation yields
$\big[P_{k+1}(T^*), |T|\big]=0$.

Replacing the pair $(U,T)$ with $(U^*,T^*)$, we get $\big[P_{k+1}(T),|T^*|\big]=0$ for $1\leq k\leq n$.
This completes the proof of \eqref{equ:property of n-centered-stronger}.
\end{proof}

\begin{lemma}\label{lem:relationship k-|T|-alpha} Let $T\in\mathcal{L}(H)$ have the polar decomposition $T=U|T|$, and let $n\in\mathbb{N}$ be such that $T$ is $(n+1)$-centered. Then for every $\alpha>0$,
\begin{equation}\label{equ:relationship k-|T|-alpha} U^k|T|^{\alpha}(U^k)^*=\big[U^k|T|(U^k)^*\big]^{\alpha}\quad \mbox{for $1\leq k\leq n+1$}.
\end{equation}
\end{lemma}

\begin{proof} Given every $k$ with $1\leq k\leq n+1$, by Lemma~\ref{lem:property of n-centered-stronger}
\begin{align*}\big[U^k|T|(U^k)^*\big]^2&=U^k|T|P_{k}(T^*)\cdot|T|(U^k)^*\\
&=\big[U^k P_{k}(T^*)\big]\cdot |T|^2\cdot (U^k)^*=U^k|T|^2(U^k)^*.
\end{align*}
The same argument shows that
\begin{align}\label{equ:n-th power of k-|T-sta|}\big[U^k|T|(U^k)^*\big]^l=U^k|T|^l(U^k)^*\quad \mbox{for all $l\in\mathbb{N}$}.
\end{align}
Let $f(t)=t^\alpha$ for $t\in [0,M]$, where $M=\big\| |T|\big\|=\|T\|$. Choose any sequence $\{P_m\}_{m=1}^{\infty}$ of polynomials such that $P_m(0)=0$ for each $m\in\mathbb{N}$, and $P_m(t)\rightarrow f(t)$ uniformly on $[0, M]$. Then from \eqref{equ:n-th power of k-|T-sta|}, we have
\begin{align*}U^k|T|^\alpha (U^k)^*
=&U^kf(|T|)(U^k)^*=\lim_{m\to\infty}U^kP_m(|T|)(U^k)^*\\
=&\lim_{m\to\infty}P_m\big[U^k|T|(U^k)^*\big]
=f\big[U^k|T|(U^k)^*\big]\\
=&\big[U^k|T|(U^k)^*\big]^\alpha.\qedhere
\end{align*}
\end{proof}

\begin{remark}\label{rem:compare with Hilbert space operators-1} In view of Lemma~\ref{lem:characterization of (n+1)-centered}, we see that  Lemma~\ref{lem:relationship k-|T|-alpha} is a generalization of \cite[Lemma~3.7~(i)]{Ito} from the Hilbert space case to the Hilbert $C^*$-module case. The method employed here is different from that in \cite{Ito}.
\end{remark}

Now, we are in the position to make a generalization of Lemma~\ref{lem:characterization of (n+1)-centered} as follows.
\begin{theorem}\label{thm:further describe of (n+2)-centered} Let $T\in\mathcal{L}(H)$ have the polar decomposition $T=U|T|$. Then for each   $n\in\mathbb{N}$, the following statements are equivalent:
\begin{enumerate}
\item[{\rm (i)}] For every $\alpha>0$ and $\beta>0$, $\big[U^k|T|^{\alpha}(U^k)^*, |T|^{\beta}\big]=0$ for $1\leq k\leq n$;
\item[{\rm (ii)}] There exist $\alpha>0$ and $\beta>0$ such that
$\big[U^k|T|^{\alpha}(U^k)^*, |T|^{\beta}\big]=0$ for $1\leq k\leq n$;
\item[{\rm (iii)}] $T$ is $(n+1)$-centered.
\end{enumerate}
\end{theorem}
\begin{proof} The implication (i)$\Longrightarrow$(ii) is obvious.

(ii)$\Longrightarrow$(iii).  Let $\alpha>0$ and $\beta>0$ be such that $\big[U^k|T|^{\alpha}(U^k)^*, |T|^{\beta}\big]=0$ for $1\leq k\leq n$.
Then for all such $k$, $\big[U^k|T|^{\alpha}(U^k)^*, |T|\big]=0$ by Lemma~\ref{lem:commutative property extended-3}~(i),  since $|T|=\big(|T|^{\beta}\big)^{\frac{1}{\beta}}$. In what follows, we prove the validity of (iii) by induction on $n$.

Case $n=1$. Suppose that $\big[U|T|^{\alpha}U^*,|T|\big]=0$. By \cite[Lemma~3.12~(i)]{Liu-Luo-Xu}
$U|T|^\alpha U^*=|T^*|^\alpha$, hence $\big[|T^*|^{\alpha},|T|\big]=0$, which leads by  Lemma~\ref{lem:commutative property extended-3}~(i) and $|T^*|=\big(|T^*|^{\alpha}\big)^{\frac{1}{\alpha}}$ to  $\big[|T^*|,|T|\big]=0$. Consequently, by \cite[Remark~4.11]{Liu-Luo-Xu-1} we conclude that $T$ is $2$-centered.

Assume that the implication (ii)$\Longrightarrow$(iii) is true for $n\in\mathbb{N}$. We show that the same is true for $n+1$. Suppose that
\begin{equation}\label{equ:2-(n+2) case of k-alpha-beta} \big[U^k|T|^{\alpha}(U^k)^*,|T|\big]=0\quad\mbox{for $1\leq k\leq n+1$}.
\end{equation}
Then obviously \eqref{equ:2-(n+2) case of k-alpha-beta} is satisfied with $n+1$ be replaced by $n$, hence $T$ is $(n+1)$-centered by the inductive hypothesis. Therefore, by Lemma~\ref{lem:relationship k-|T|-alpha}, \eqref{equ:2-(n+2) case of k-alpha-beta} can be rewritten as
\begin{equation*} \big[(A_k)^\alpha,|T|\big]=0\quad\mbox{for $1\leq k\leq n+1$},
\end{equation*}
in which $A_k$ is defined by
\begin{equation}\label{equ:defn of A k}A_k=U^k|T|(U^k)^*.\end{equation}
Utilizing Lemma~\ref{lem:commutative property extended-3}~(i) we see that
\begin{equation*}\big[A_k,|T|\big]=0\quad\mbox{for $1\leq k\leq n+1$},
\end{equation*}
which means by Lemma~\ref{lem:characterization of (n+1)-centered} that $T$ is $(n+2)$-centered, as desired.

(iii)$\Longrightarrow$(i). Let $A_k$ be defined by \eqref{equ:defn of A k}.
By Lemmas~\ref{lem:characterization of (n+1)-centered} and \ref{lem:relationship k-|T|-alpha},
$$\big[A_k,|T|\big]=0\quad\mbox{and}\quad (A_k)^\alpha=U^k|T|^\alpha(U^k)^*$$ for each $k=1,2,\cdots, n$ and every $\alpha>0$.
Utilizing the first equation above and Lemma~\ref{lem:commutative property extended-3}~(i) we see that for every $\alpha>0$ and $\beta>0$,
$$\big[U^k|T|^\alpha(U^k)^*, |T|^{\beta}\big]=\big[(A_k)^\alpha, |T|^\beta\big]=0.\qedhere$$
\end{proof}

\begin{remark}\label{rem:one result of (n+1)-centered}From item (i) of the preceding theorem and Lemma~\ref{lem:commutative property extended-3}~(ii), it can be concluded that
\begin{equation}\label{equ:one result of (n+1)-centered-alpha}\big[U^k|T|^\alpha(U^k)^*, U^*U\big]=0\quad\mbox{for every $\alpha>0$ and $1\leq k\leq n$},\end{equation} since by Lemma~\ref{lem:Range Closure of Ta and T} we have $\overline{\mathcal{R}(|T|^\beta)}=\overline{\mathcal{R}(|T|)}=\mathcal{R}(U^*U)$.
\end{remark}

\section{Generalized aluthge transforms and $(n+1)$-centered operators}\label{sec:generalized Aluthge transform}

Let $T\in\mathcal{L}(H)$ have the polar decomposition $T=U|T|$, and let $\alpha$ and $\beta$ be positive numbers. The associated operator $\widetilde{T}_{\alpha, \beta}$ defined by
\begin{equation}\label{equ:defn of generalized Aluthge transforms} \widetilde{T}_{\alpha, \beta}=|T|^{\alpha}U|T|^{\beta}
\end{equation}
is called the generalized Aluthge transform of $T$. Specifically, $\widetilde{T}_{\frac12,\frac12}$ is usually denoted by $\widetilde{T}$ and is called
the Aluthge trasformation of $T$ \cite{Aluthge}.

Recall that an operator $T\in\mathcal{L}(H)$ is said to be binormal \cite{Campbell-1} if $T^*T$ and $TT^*$ are commutative. That is, $[T^*T,TT^*]=0$.
When $T$ has the polar decomposition $T=U|T|$, it is known that $T$ is binormal if and only if $T$ is $2$-centered (see e.g.,\cite[Theorem~3.1]{Ito-Yamazaki-Yanagida} and \cite[Remark~4.11]{Liu-Luo-Xu-1}).

The purpose of this section is to give additional characterizations of $(n+1)$-centered operators in terms of generalized Aluthge transforms; see Theorems~\ref{thm:cited lemma of Ito++} and \ref{thm:characterization of (n+1)-centered when binormal}  for the details.

\begin{lemma}\label{lem:T-alpha-beta separate into T-beta and T-alpha} Let $T\in\mathcal{L}(H)$ have the polar decomposition $T=U|T|$, and let $n\in\mathbb{N}$ be such that $T$ is $(n+1)$-centered. Then for every $\alpha>0$ and $\beta>0$, and every $k\in\mathbb{N}$ with $1\leq k\leq n+2$,
\begin{equation}\label{equ:T-alpha-beta separate into T-beta and T-alpha} U^k|\widetilde{T}_{\alpha, \beta}|(U^k)^*=U^k|T|^{\beta}(U^k)^*\cdot U^{k-1}|T|^{\alpha}(U^{k-1})^*,
\end{equation}
where $\widetilde{T}_{\alpha,\beta}$ is defined by \eqref{equ:defn of generalized Aluthge transforms}.
\end{lemma}
\begin{proof} Let $\alpha$ and $\beta$ be any positive numbers, and let $P_k(T)$ and $P_k(T^*)$ be defined by \eqref{equ:defn of P n T T-sta} for every $k\in\mathbb{N}$. Suppose that $2\leq k\leq n+2$. Then $1\leq k-1\leq n+1$, so the first equation in \eqref{equ:property of n-centered-stronger} together with Lemma~\ref{lem:commutative property extended-3}~(i) yields
\begin{equation}\label{equ:half-property of n-centered-stronger}\big[P_{k-1}(T^*),|T|^{\alpha}\big]=0.
\end{equation}
Since $T$ is $(n+1)$-centered, it is $2$-centered, hence is binormal, so by \cite[Theorem~4.14]{Liu-Luo-Xu-1} we have
\begin{align}\label{eqn:formula for square root of widetilde T-alpha-beta}|&\widetilde{T}_{\alpha,\beta}|=U^*|T|^{\alpha} U\cdot |T|^{\beta}=|T|^{\beta}\cdot U^*|T|^{\alpha} U,\\
\label{eqn:formula for square root of widetilde T-alpha-beta-star}&|\widetilde{T}_{\alpha,\beta}^*|=|T|^{\alpha}\cdot |T^*|^{\beta}=|T^*|^{\beta}\cdot |T|^{\alpha},\\
\label{equ:basic commutativity derived from binormal condition}&\big[|T|^{\alpha},UU^*\big]=\big[U^*U,|T^*|^{\beta}\big]=\big[U^*U,UU^*\big]=0,\\
\label{equ:commutative related binormal-1}&\big[U|T|^\alpha U^*,|T|^\beta\big]=\big[U^*|T|^{\alpha}U,|T|^{\beta}\big]=0.
\end{align}
Therefore, for each $k\in\mathbb{N}$ with $2\leq k\leq n+2$ we have
\begin{align*}U^k|\widetilde{T}_{\alpha,\beta}|(U^k)^*
&=U^k\big(|T|^{\beta}U^*|T|^{\alpha}U\big)(U^k)^*\quad \mbox{\big(by \eqref{eqn:formula for square root of widetilde T-alpha-beta}\big)}\\
&=U^k|T|^{\beta}U^*\cdot|T|^{\alpha}UU^*\cdot(U^{k-1})^*\\
&=U^k|T|^{\beta}U^*\cdot UU^*|T|^{\alpha}(U^{k-1})^*\quad \mbox{\big(by \eqref{equ:basic commutativity derived from binormal condition}\big)}\\
&=U^k|T|^{\beta}U^*|T|^{\alpha}P_{k-1}(T^*)\cdot(U^{k-1})^*\\
&=U^k|T|^{\beta}U^*P_{k-1}(T^*)|T|^{\alpha}(U^{k-1})^*\quad \mbox{\big(by \eqref{equ:half-property of n-centered-stronger}\big)}\\
&=U^k|T|^{\beta}(U^k)^*\cdot U^{k-1}|T|^{\alpha}(U^{k-1})^*.
\end{align*}

In the remaining case of $k=1$, by \eqref{eqn:formula for square root of widetilde T-alpha-beta} and \eqref{equ:basic commutativity derived from binormal condition} we obtain
\begin{align*}U|\widetilde{T}_{\alpha,\beta}|U^*&=U\big(|T|^{\beta}U^*|T|^{\alpha}U\big)U^*
=U|T|^{\beta}U^*UU^*|T|^{\alpha}=U|T|^{\beta}U^*|T|^{\alpha}.
\end{align*}
This completes the proof of \eqref{equ:T-alpha-beta separate into T-beta and T-alpha}.
\end{proof}

\begin{remark}\label{rem:compare with Hilbert space operators-3} The special case of the preceding lemma is considered in \cite[Lemma~3.7~(vii)]{Ito} for Hilbert space operators, where $k=n+1$ and $\alpha=\beta=\frac12$.
\end{remark}

Now, we provide the technical result of this section as follows.

\begin{theorem}\label{thm:cited lemma of Ito++} Suppose that $T$ is binormal and has the polar decomposition $T=U|T|$. Let $\widetilde{T}_{\alpha,\beta}$ be defined by \eqref{equ:defn of generalized Aluthge transforms} for $\alpha>0$ and $\beta>0$. Then for every $n\in\mathbb{N}$ with $n\ge 2$, the following statements are equivalent:
\begin{enumerate}
\item[{\rm (i)}] $T$ is $(n+1)$-centered;
\item[{\rm (ii)}] $\big[U^k|T|(U^k)^*,|T|\big]=0$ for $2\leq k\leq n$;
\item[{\rm (iii)}] For every $\alpha>0$ and $\beta>0$,
$$\big[U^k|\widetilde{T}_{\alpha,\beta}|(U^k)^*, |\widetilde{T}_{\alpha,\beta}|\big]=0\quad \mbox{for $1\leq k\leq n-1$;}$$
\item[{\rm (iv)}] There exist $\alpha>0$ and $\beta>0$ such that
\begin{equation*}\big[U^k|\widetilde{T}_{\alpha,\beta}|(U^k)^*, |\widetilde{T}_{\alpha,\beta}|\big]=0\quad \mbox{for $1\leq k\leq n-1$}.
\end{equation*}
\end{enumerate}
\end{theorem}
\begin{proof} Since $T$ is binormal, by \cite[Theorem~4.14]{Liu-Luo-Xu-1}  $|\widetilde{T}_{\alpha, \beta}|$ is given by \eqref{eqn:formula for square root of widetilde T-alpha-beta}.

 (i)$\Longleftrightarrow$(ii).  Note that $U|T|U^*=|T^*|$ \cite[Lemma~3.12]{Liu-Luo-Xu}, so by the binormality of $T$ and Lemma~\ref{lem:commutative property extended-3}~(i) we have
 $$\big[U|T|U^*,|T|\big]=\big[|T^*|,|T|\big]=0,$$
which leads by Lemma~\ref{lem:characterization of (n+1)-centered} to the equivalence of (i) and (ii).

(i)+(ii)$\Longrightarrow$(iii). Let
\begin{equation}\label{equ:defn of A k and B k} A_k=U^k|\widetilde{T}_{\alpha,\beta}|(U^k)^*\cdot|\widetilde{T}_{\alpha,\beta}|\quad \mbox{and}\quad B_k=|\widetilde{T}_{\alpha,\beta}|\cdot U^k|\widetilde{T}_{\alpha,\beta}|(U^k)^*
\end{equation}
for $1\le k\le n-1$.
As  $k+1=2,3,\cdots,n$, we can use Theorem~\ref{thm:further describe of (n+2)-centered}~(i) and \eqref{equ:one result of (n+1)-centered-alpha} repeatedly to derive equations (except for the equations specified) as
  \begin{align*} A_k&=U^k|T|^{\beta}(U^k)^*\cdot U^{k-1}|T|^{\alpha}(U^{k-1})^*\cdot|T|^{\beta}\cdot U^*|T|^{\alpha}U\quad\mbox{\big(by \eqref{equ:T-alpha-beta separate into T-beta and T-alpha} and \eqref{eqn:formula for square root of widetilde T-alpha-beta}\big)} \\
&=|T|^{\beta}\cdot U^k|T|^{\beta}(U^k)^*\cdot U^{k-1}|T|^{\alpha}(U^{k-1})^*\cdot U^*|T|^{\alpha}U\\
&=|T|^{\beta}\cdot U^k|T|^{\beta}(U^k)^*\cdot U^{k-1}|T|^{\alpha}(U^{k-1})^*\cdot (U^*U)^2U^*|T|^{\alpha}U\\
&=|T|^{\beta}\cdot (U^*U) U^k|T|^{\beta}(U^k)^*\cdot (U^*U) U^{k-1}|T|^{\alpha}(U^{k-1})^*\cdot U^*|T|^{\alpha}U\\
&=|T|^{\beta}U^*\cdot U^{k+1}|T|^{\beta}(U^{k+1})^*\cdot U^k|T|^{\alpha}(U^k)^*\cdot|T|^{\alpha}U\\
&=|T|^{\beta}U^*|T|^{\alpha}\cdot U^{k+1}|T|^{\beta}(U^{k+1})^*\cdot U^k|T|^{\alpha}(U^k)^*U\\
 &=|\widetilde{T}_{\alpha,\beta}|\cdot U^k|T|^{\beta}(U^k)^*\cdot U^*U\cdot U^{k-1}|T|^{\alpha}(U^{k-1})^*\cdot U^*U\quad \mbox{\big(by \eqref{eqn:formula for square root of widetilde T-alpha-beta}\big)}\\
&=\left[|\widetilde{T}_{\alpha,\beta}|(U^*U)^2\right]\cdot U^k|T|^{\beta}(U^k)^*\cdot U^{k-1}|T|^{\alpha}(U^{k-1})^*\\
&=|\widetilde{T}_{\alpha,\beta}|\cdot U^k|T|^{\beta}(U^k)^*\cdot U^{k-1}|T|^{\alpha}(U^{k-1})^*\quad \mbox{\big(by \eqref{eqn:formula for square root of widetilde T-alpha-beta}\big)}\\
&=|\widetilde{T}_{\alpha,\beta}|\cdot U^k|\widetilde{T}_{\alpha,\beta}|(U^k)^*\quad\mbox{\big(by \eqref{equ:T-alpha-beta separate into T-beta and T-alpha}\big)}\\
&=B_k.
\end{align*}

The implication (iii)$\Longrightarrow$(iv) is obvious.

(iv)$\Longrightarrow$(ii). Let $\alpha>0$ and $\beta>0$ be such that item~(iv) of this theorem is satisfied. We prove the validity of (ii) by induction on $n$.

Case $n=2$. Suppose that $\big[U|\widetilde{T}_{\alpha,\beta}|U^*,|\widetilde{T}_{\alpha,\beta}|\big]=0$. Then clearly $A_1=B_1$, where $A_1$ and $B_1$ are  defined by \eqref{equ:defn of A k and B k} for $k=1$.
Since $T$ is binormal, we may combine \eqref{eqn:formula for square root of widetilde T-alpha-beta} with \eqref{equ:basic commutativity derived from binormal condition} and \cite[Lemma~3.12~(i)]{Liu-Luo-Xu} to get
\begin{align*} U|\widetilde{T}_{\alpha,\beta}|U^*&=U\big(|T|^{\beta}\cdot U^*|T|^{\alpha}U\big)U^*
=U|T|^{\beta}\cdot U^*UU^*|T|^{\alpha}\\
&=U|T|^{\beta}U^*\cdot|T|^{\alpha}=|T^*|^{\beta}\cdot|T|^{\alpha}.
\end{align*}
The equations above, together with  \eqref{eqn:formula for square root of widetilde T-alpha-beta}, \eqref{equ:basic commutativity derived from binormal condition} and \eqref{equ:commutative related binormal-1} yield
\begin{align*}A_1&=U|\widetilde{T}_{\alpha,\beta}|U^*\cdot|\widetilde{T}_{\alpha,\beta}|=U|T|^{\beta}U^*\cdot|T|^{\alpha}\cdot (U^*U)^2U^*|T|^{\alpha}U|T|^{\beta}\\
&=U^*U\cdot U|T|^{\beta}U^*\cdot U^*U|T|^{\alpha}\cdot U^*|T|^{\alpha}U|T|^{\beta}\\
&=U^*\cdot U^2|T|^{\beta}(U^2)^*\cdot U|T|^{\alpha}U^*\cdot|T|^{\alpha}U|T|^{\beta}\\
&=U^*\cdot U^2|T|^{\beta}(U^2)^*\cdot|T|^{\alpha}\cdot U|T|^{\alpha}U^*\cdot U|T|^{\beta}\\
&=U^*\cdot U^2|T|^{\beta}(U^2)^*|T|^{\alpha}\cdot U|T|^{\alpha}|T|^{\beta}.
\end{align*}
Similarly,
\begin{align*}B_1&=|\widetilde{T}_{\alpha,\beta}|\cdot U|\widetilde{T}_{\alpha,\beta}|U^*=U^*|T|^{\alpha}U\cdot|T|^{\beta}\cdot U|T|^{\beta}U^*\cdot|T|^{\alpha}\\
&=U^*|T|^{\alpha}U\cdot U|T|^{\beta}U^*\cdot|T|^{\alpha}\cdot|T|^{\beta}\\
&=U^*\cdot|T|^{\alpha}U^2|T|^{\beta}(U^2)^*\cdot U|T|^{\alpha}|T|^{\beta}.
\end{align*}
The expressions of $A_1$ and $B_1$ as above, together with $UA_1=UB_1$ and \eqref{equ:basic commutativity derived from binormal condition}, yield
\begin{align*}U^2|T|^{\beta}(U^2)^*|T|^{\alpha}\cdot U|T|^{\alpha}|T|^{\beta}
=|T|^{\alpha}U^2|T|^{\beta}(U^2)^*\cdot U|T|^{\alpha}|T|^{\beta},
\end{align*}
which menas that
\begin{align*}U^2|T|^{\beta}(U^2)^*|T|^{\alpha}\xi
=|T|^{\alpha}U^2|T|^{\beta}(U^2)^*\xi
\end{align*}
for every $\xi\in\overline{\mathcal{R}\big(U|T|^{\alpha+\beta}\big)}$. It follows that
\begin{equation*} U^2|T|^{\beta}(U^2)^*|T|^{\alpha}\cdot UU^*=|T|^{\alpha}U^2|T|^{\beta}(U^2)^*\cdot UU^*,
\end{equation*}
since by  Lemmas~\ref{lem:rang characterization-1} and  \ref{lem:Range Closure of Ta and T} we have $\overline{\mathcal{R}(U|T|^{\alpha+\beta})}=\overline{\mathcal{R}(UU^*)}=\mathcal{R}(UU^*)$.
Hence, by \eqref{equ:relationship k-|T|-alpha}, \eqref{equ:basic commutativity derived from binormal condition} and \eqref{equ:commutative related binormal-1} we see that
\begin{align*}\big[U^2|T|(U^2)^*\big]^{\beta}\cdot|T|^{\alpha}&=U^2|T|^{\beta}(U^2)^*\cdot|T|^{\alpha}=U^2|T|^{\beta}(U^2)^*\cdot UU^*|T|^{\alpha}\\
&=U^2|T|^{\beta}(U^2)^*|T|^{\alpha}\cdot UU^*=|T|^{\alpha}U^2|T|^{\beta}(U^2)^*\cdot UU^*\\
&=|T|^{\alpha}\cdot U^2|T|^{\beta}(U^2)^*=|T|^{\alpha}\cdot \big[U^2|T|(U^2)^*\big]^{\beta}.
\end{align*}
Thus, $\big[U^2|T|(U^2)^*,|T|\big]=0$ by Lemma~\ref{lem:commutative property extended-3}~(i). This completes the proof of the case $n=2$.

Assume that the implication (iv)$\Longrightarrow$(ii) is true for $n\in\mathbb{N}$. We show that the same is true for $n+1$. Suppose that
\begin{equation}\label{equ:commutativity in terms of U-k-alpha-beta-11}\big[U^k|\widetilde{T}_{\alpha,\beta}|(U^k)^*,|\widetilde{T}_{\alpha,\beta}|\big]=0\quad\mbox{for $1\leq k\leq n$}.
\end{equation}
Then obviously \eqref{equ:commutativity in terms of U-k-alpha-beta-11} is satisfied with $n$ be replaced by $n-1$, hence by the inductive hypothesis one has
$$\big[U^k|T|(U^k)^*,|T|\big]=0\quad\mbox{for $2\leq k\leq n$}.$$
So, $T$ is $(n+1)$-centered by the equivalence of (i) and (ii). Thus, we need only to prove that
\begin{equation}\label{equ:commutativity U-n+2-|T|-111}\big[U^{n+1}|T|(U^{n+1})^*,|T|\big]=0.
\end{equation}

To this end, we consider the operators $A_n$ and $B_n$ defined by  \eqref{equ:defn of A k and B k} for $k=n$ therein. Since $T$ is $(n+1)$-centered, by \eqref{equ:T-alpha-beta separate into T-beta and T-alpha}, \eqref{eqn:formula for square root of widetilde T-alpha-beta}, \eqref{equ:one result of (n+1)-centered-alpha} and Theorem~\ref{thm:further describe of (n+2)-centered}~(i) we have
\begin{align*}A_n&=U^n|T|^{\beta}(U^n)^*\cdot U^{n-1}|T|^{\alpha}(U^{n-1})^*\cdot(U^*U)^2U^*|T|^{\alpha}U|T|^{\beta}\\
&=U^*U\cdot U^n|T|^{\beta}(U^n)^*\cdot U^*U\cdot U^{n-1}|T|^{\alpha}(U^{n-1})^*\cdot U^*|T|^{\alpha}U|T|^{\beta}\\
&=U^*\cdot U^{n+1}|T|^{\beta}(U^{n+1})^*\cdot U^n|T|^{\alpha}(U^n)^*\cdot |T|^{\alpha}U|T|^{\beta}\\
&=U^*\cdot U^{n+1}|T|^{\beta}(U^{n+1})^*\cdot|T|^{\alpha}\cdot U^n|T|^{\alpha}(U^n)^*\cdot U|T|^{\beta}\\
&=U^*\cdot U^{n+1}|T|^{\beta}(U^{n+1})^*|T|^{\alpha}\cdot U^n|T|^{\alpha}(U^{n-1})^*|T|^{\beta}.
\end{align*}
Similarly,
\begin{align*}B_n&=U^*|T|^{\alpha}U|T|^{\beta}\cdot U^n|T|^{\beta}(U^n)^*\cdot U^{n-1}|T|^{\alpha}(U^{n-1})^*\\
&=U^*|T|^{\alpha}U(U^*U)\cdot U^n|T|^{\beta}(U^n)^*\cdot U^{n-1}|T|^{\alpha}(U^{n-1})^*\cdot |T|^{\beta}\\
&=U^*|T|^{\alpha}U\cdot U^n|T|^{\beta}(U^n)^*\cdot U^*U\cdot U^{n-1}|T|^{\alpha}(U^{n-1})^*\cdot|T|^{\beta}\\
&=U^*\cdot|T|^{\alpha}U^{n+1}|T|^{\beta}(U^{n+1})^*\cdot U^n|T|^{\alpha}(U^{n-1})^*|T|^{\beta}.
\end{align*}
Utilizing \eqref{equ:commutativity in terms of U-k-alpha-beta-11}  in the case of $k=n$, we obtain $A_n=B_n$.
It follows from $UA_n=UB_n$ and \eqref{equ:basic commutativity derived from binormal condition} that
\begin{equation}\label{equ:commutativity with Cn} U^{n+1}|T|^{\beta}(U^{n+1})^*|T|^{\alpha}\cdot C_n=|T|^{\alpha}U^{n+1}|T|^{\beta}(U^{n+1})^*\cdot C_n,
\end{equation}
where $C_n=U^n|T|^{\alpha}(U^{n-1})^*|T|^{\beta}$.
Since $U^*U|T|^{\beta}=|T|^{\beta}$, by Lemmas~\ref{lem:rang characterization-1} and \ref{lem:Range Closure of Ta and T} we have
\begin{align*}\overline{\mathcal{R}(C_n)}
&=\overline{\mathcal{R}\big[U^n|T|^{\alpha}(U^n)^*U|T|^{\beta}\big]}
=\overline{\mathcal{R}\big[U^n|T|^{\alpha}(U^n)^*UU^*\big]}\\
&=\overline{\mathcal{R}\big[U^n|T|^{\alpha}(U^n)^*\big]}
=\overline{\mathcal{R}\big[U^n|T|^{\frac{\alpha}{2}}\big]}
=\overline{\mathcal{R}(U^nU^*)}.
\end{align*}
This, together with \eqref{equ:commutativity with Cn}, yields
\begin{equation}\label{equ:commutativity with UnU*}U^{n+1}|T|^{\beta}(U^{n+1})^*|T|^{\alpha}U^nU^*=|T|^{\alpha}U^{n+1}|T|^{\beta}(U^{n+1})^* U^nU^*.
\end{equation}
Let $P_n(T)$ be defined by \eqref{equ:defn of P n T T-sta}. Since $T$ is $(n+1)$-centered, by \cite[Lemma~5.6]{Liu-Luo-Xu-1}  $\big[P_n(T),|T|\big]=0$, which gives $[P_n(T), |T|^{\alpha}]=0$ by Lemma~\ref{lem:commutative property extended-3}~(i). Hence,
\begin{align*}U^{n+1}|T|^{\beta}(U^{n+1})^*\cdot|T|^{\alpha}
&=U^{n+1}|T|^{\beta}U^*\cdot(U^{n})^*P_{n}(T)|T|^{\alpha}\\
&=U^{n+1}|T|^{\beta}(U^{n+1})^*\cdot|T|^{\alpha}P_{n}(T)\\
&=U^{n+1}|T|^{\beta}(U^{n+1})^*|T|^{\alpha}U^nU^*\cdot(U^{n-1})^*\\
&=|T|^{\alpha}U^{n+1}|T|^{\beta}(U^{n+1})^*U^nU^*\cdot(U^{n-1})^*\ \mbox{\big(by \eqref{equ:commutativity with UnU*}\big)}\\
&=|T|^{\alpha}\cdot U^{n+1}|T|^{\beta}(U^{n+1})^*.
\end{align*}
So, from \eqref{equ:relationship k-|T|-alpha} we obtain
$$\left[U^{n+1}|T|(U^{n+1})^*\right]^{\beta}\cdot|T|^{\alpha}=|T|^{\alpha}\left[U^{n+1}|T|(U^{n+1})^*\right]^{\beta},$$
which clearly leads by Lemma~\ref{lem:commutative property extended-3}~(i) to \eqref{equ:commutativity U-n+2-|T|-111}.
\end{proof}

Suppose that $T\in\mathcal{L}(H)$ has the polar decomposition $T=U|T|$, let

\begin{equation}\label{equ:defn of widetilde U} \widetilde{U}=U^*U^2.
\end{equation}
A simple calculation shows that
\begin{align}\label{equ:formulas for n-th power of widetilde U}\widetilde{U}^n=U^*U^{n+1}\quad\mbox{for every $n\in\mathbb{N}$}.
\end{align}

\begin{lemma} \label{lem:delete U star U} Let $T\in\mathcal{L}(H)$ have the polar decomposition $T=U|T|$, and let $n\in\mathbb{N}$ be such that $T$ is $(n+1)$-centered. Then for every $\alpha>0$ and $\beta>0$,
\begin{align}\label{eqn:important relationship--1}
\widetilde{U}^k|\widetilde{T}_{\alpha,\beta}|(\widetilde{U}^k)^*\cdot |\widetilde{T}_{\alpha,\beta}|&=U^k|\widetilde{T}_{\alpha,\beta}|(U^k)^*\cdot |\widetilde{T}_{\alpha,\beta}|\quad \mbox{for $1\leq k\leq n$},
\end{align}
where $\widetilde{T}_{\alpha,\beta}$ and $\widetilde{U}$ are defined by \eqref{equ:defn of generalized Aluthge transforms} and \eqref{equ:defn of widetilde U}, respectively.
\end{lemma}
\begin{proof} By Lemma~\ref{lem:Range Closure of Ta and T} $\overline{\mathcal{R}(|T|^\alpha)}=\overline{\mathcal{R}(|T|)}=\overline{\mathcal{R}(T^*)}$, which means that
$$U^*U|T|^\alpha=|T|^\alpha=|T|^\alpha U^*U.$$
Thus, $\big[|T|^\alpha,U^*U\big]=0$. Furthermore, from \eqref{equ:one result of (n+1)-centered-alpha} we have
$$\big[U^k|T|^\alpha(U^k)^*, U^*U\big]=0\quad\mbox{and}\quad \big[U^k|T|^\beta(U^k)^*, U^*U\big]=0$$
for $1\leq k\leq n$. This shows that  for $1\leq k\leq n$,
$$\big[U^k|T|^{\beta}(U^k)^*\cdot U^{k-1}|T|^{\alpha}(U^{k-1})^*,U^*U\big]=0.$$
Consequently,  by \eqref{equ:T-alpha-beta separate into T-beta and T-alpha} we obtain
\begin{equation}\label{equ:commutativity in terms of U-star U-22}\big[U^k|\widetilde{T}_{\alpha,\beta}| (U^k)^*, U^*U\big]=0\quad \mbox{for $1\le k\le n$}.\end{equation}
Hence, for  $k=1,2,\cdots, n$,
\begin{align*}\widetilde{U}^k|\widetilde{T}_{\alpha,\beta}|(\widetilde{U}^k)^*\cdot |\widetilde{T}_{\alpha,\beta}|&=(U^*U)U^k|\widetilde{T}_{\alpha,\beta}| (U^k)^*\cdot (U^*U)|\widetilde{T}_{\alpha,\beta}|\quad \mbox{\big(by \eqref{equ:formulas for n-th power of widetilde U}\big)}\\
&=U^k|\widetilde{T}_{\alpha,\beta}| (U^k)^*\cdot (U^*U)(U^*U)|\widetilde{T}_{\alpha,\beta}|\quad \mbox{\big(by \eqref{equ:commutativity in terms of U-star U-22}\big)}\\
&=U^k|\widetilde{T}_{\alpha,\beta}| (U^k)^*\cdot |\widetilde{T}_{\alpha,\beta}|\quad \mbox{\big(by \eqref{eqn:formula for square root of widetilde T-alpha-beta}\big).}
\end{align*}
This completes the proof of \eqref{eqn:important relationship--1}.
\end{proof}

We provide an additional characterization of the $(n+1)$-centered operator as follows.
\begin{theorem}\label{thm:characterization of (n+1)-centered when binormal} Let $T\in\mathcal{L}(H)$ have the polar decomposition $T=U|T|$, and let $\widetilde{T}_{\alpha,\beta}$ be defined by \eqref{equ:defn of generalized Aluthge transforms} for every $\alpha>0$ and $\beta>0$. Then for each  $n\in\mathbb{N}$, the following statements are equivalent:
\begin{enumerate}
\item[{\rm (i)}] $T$ is $(n+1)$-centered;
\item[{\rm (ii)}] $T$ is binormal and for every $\alpha>0, \beta>0$, $\widetilde{T}_{\alpha,\beta}$ is $n$-centered;
\item[{\rm (iii)}] $T$ is binormal and there exist $\alpha>0$ and $\beta>0$ such that $\widetilde{T}_{\alpha,\beta}$ is $n$-centered.
\end{enumerate}
\end{theorem}
\begin{proof} (i)$\Longrightarrow$(ii). Assume that $T$ is $(n+1)$-centered. Then obviously $T$ is $2$-centered, hence $T$ is binormal. So for every $\alpha>0$ and $\beta>0$, by \cite[Theorem~4.14]{Liu-Luo-Xu-1}  $\widetilde{T}_{\alpha,\beta}=\widetilde{U}|\widetilde{T}_{\alpha,\beta}|$ is the polar decomposition such that $|\widetilde{T}_{\alpha, \beta}|$ is given by \eqref{eqn:formula for square root of widetilde T-alpha-beta},
where $\widetilde{U}$ is defined by \eqref{equ:defn of widetilde U}. This shows the validity of the implication (i)$\Longrightarrow$(ii) in the case that $n=1$.

Suppose that $n\ge 2$. By Theorem~\ref{thm:cited lemma of Ito++}~(iii), we have
\begin{equation*}\big[U^k|\widetilde{T}_{\alpha,\beta}|(U^k)^*,|\widetilde{T}_{\alpha,\beta}|\big]=0\quad \mbox{for $1\leq k\leq n-1$},
\end{equation*}
which is combined with Lemma~\ref{lem:delete U star U} to get
\begin{equation*}\big[\widetilde{U}^k|\widetilde{T}_{\alpha,\beta}|(\widetilde{U}^k)^*,|\widetilde{T}_{\alpha,\beta}|\big]=0\quad \mbox{for $1\leq k\leq n-1$}.
\end{equation*}
Applying Lemma~\ref{lem:characterization of (n+1)-centered} to the polar decomposition $\widetilde{T}_{\alpha,\beta}=\widetilde{U}|\widetilde{T}_{\alpha,\beta}|$, we conclude that $\widetilde{T}_{\alpha,\beta}$ is $n$-centered.

The implication (ii)$\Longrightarrow$(iii) is obvious.

(iii)$\Longrightarrow$(i). Suppose that  $T$ is binormal. Then $T$ is $2$-centered, hence the implication (iii)$\Longrightarrow$(i)
 is valid in the case of $n=1$. Since $T$ is binormal, $\widetilde{T}_{\alpha,\beta}=\widetilde{U}|\widetilde{T}_{\alpha,\beta}|$ is the polar decomposition. In what follows, we prove by induction on $n$ that $T$ will be $(n+1)$-centered whenever $\widetilde{T}_{\alpha,\beta}$ is $n$-centered for some $\alpha>0$ and $\beta>0$.

The case $n=1$ has already been proved. Assume the assertion is true when  $\widetilde{T}_{\alpha,\beta}$ is $n$-centered.
Suppose now that $\widetilde{T}_{\alpha,\beta}$ is $(n+1)$-centered. Then by Lemma~\ref{lem:characterization of (n+1)-centered} we have
\begin{equation}\label{equ:T-alpha-beta is (n+1)-centered}\big[\widetilde{U}^k|\widetilde{T}_{\alpha, \beta}|(\widetilde{U}^k)^*,|\widetilde{T}_{\alpha,\beta}|\big]=0\quad \mbox{for $1\leq k\leq n$}.
\end{equation}
Also, we know by the inductive hypothesis that $T$ is $(n+1)$-centered, since $\widetilde{T}_{\alpha,\beta}$ is obviously $n$-centered. We can therefore use \eqref{equ:T-alpha-beta is (n+1)-centered} and Lemma~\ref{lem:delete U star U} to obtain
\begin{equation*}\big[U^k|\widetilde{T}_{\alpha,\beta}|(U^k)^*,|\widetilde{T}_{\alpha,\beta}|\big]=0\quad \mbox{for $1\leq k\leq n$},
\end{equation*}
which is combined with Theorem~\ref{thm:cited lemma of Ito++} to get
\begin{equation*}\big[U^k|T|(U^k)^*,|T|\big]=0\quad \mbox{for $1\leq k\leq n+1$}.
\end{equation*}
Therefore, by Lemma~\ref{lem:characterization of (n+1)-centered} $T$ is $(n+2)$-centered. This completes the proof of the implication (iii)$\Longrightarrow$(i).
\end{proof}

\section{Generalized iterated aluthge transforms and $(n+1)$-centered operators}\label{sec:iterated Aluthge transforms}

Suppose that $T\in\mathcal{L}(H)$ has the polar decomposition $T=U|T|$. Let $\widetilde{T}_{\alpha,\beta}$ be defined by \eqref{equ:defn of generalized Aluthge transforms} for $\alpha>0$ and $\beta>0$. The $n$-th generalized Aluthge transform of $T$ is defined inductively as $\widetilde{T}_{\alpha,\beta}^{(0)}=T$, $\widetilde{T}_{\alpha,\beta}^{(1)}=\widetilde{T}_{\alpha,\beta}$, and for each $n\in\mathbb{N}$ with $n\ge 2$, $\widetilde{T}_{\alpha,\beta}^{(n)}=\widetilde{\big(\widetilde{T}_{\alpha,\beta}^{(n-1)}\big)}_{\alpha,\beta}$ whenever $\widetilde{T}_{\alpha,\beta}^{(n-1)}$ has the polar decomposition.

The purpose of this section is to give characterizations of $(n+1)$-centered operators in terms of generalized iterated Aluthge transforms.

\begin{theorem}\label{thm:binormality of T-alpha-beta-n}
 Let $T\in\mathcal{L}(H)$ have the polar decomposition $T=U|T|$, and $n\in\mathbb{N}$ with $n\ge 2$.  Then the following statements are equivalent:
\begin{enumerate}
\item[{\rm (i)}] $T$ is $(n+1)$-centered;
\item[{\rm (ii)}] For every $\alpha>0$ and $\beta>0$, $\widetilde{T}_{\alpha,\beta}^{(k)}$ is $2$-centered for  $0\le k\le n-1$;
\item[{\rm (iii)}] There exist $\alpha>0$ and $\beta>0$ such that $\widetilde{T}_{\alpha,\beta}^{(k)}$ is $2$-centered for $0\le k\le n-1$.
\end{enumerate}
\end{theorem}
\begin{proof} (i)$\Longrightarrow$(ii). A repeatedly using the implication (i)$\Longrightarrow$(ii) of Theorem~\ref{thm:characterization of (n+1)-centered when binormal} yields
\begin{align*}
\mbox{$T$ is $(n+1)$-centered}&\Longrightarrow \mbox{$T$ is binormal, $\widetilde{T}_{\alpha,\beta}$ is $n$-centered},\\
\mbox{$\widetilde{T}_{\alpha,\beta}$ is $n$-centered}&\Longrightarrow \mbox{$\widetilde{T}_{\alpha,\beta}$ is binormal, $\widetilde{T}_{\alpha,\beta}^{(2)}$ is $(n-1)$-centered},\\
\mbox{$\widetilde{T}_{\alpha,\beta}^{(2)}$ is $(n-1)$-centered}&\Longrightarrow \mbox{$\widetilde{T}_{\alpha,\beta}^{(2)}$ is binormal, $\widetilde{T}_{\alpha,\beta}^{(3)}$ is $(n-2)$-centered},\\
&\quad \vdots\\
\mbox{$\widetilde{T}_{\alpha,\beta}^{(n-2)}$ is $3$-centered}&\Longrightarrow \mbox{$\widetilde{T}_{\alpha,\beta}^{(n-2)}$ is binormal, $\widetilde{T}_{\alpha,\beta}^{(n-1)}$ is $2$-centered},\\
\mbox{$\widetilde{T}_{\alpha,\beta}^{(n-1)}$ is $2$-centered}&\Longrightarrow \mbox{$\widetilde{T}_{\alpha,\beta}^{(n-1)}$ is binormal and $\widetilde{T}_{\alpha,\beta}^{(n)}$ is meaningful,}
\end{align*}
since $\widetilde{T}_{\alpha,\beta}^{(n)}$ is meaningful if and only if $\widetilde{T}_{\alpha,\beta}^{(n-1)}$ has the polar decomposition.

The implication (ii)$\Longrightarrow$(iii) is clear.

(iii)$\Longrightarrow$(i). By assumption $\widetilde{T}_{\alpha,\beta}^{(n-1)}$ has the polar decomposition, so a repeatedly using the implication (iii)$\Longrightarrow$(i) of Theorem~\ref{thm:characterization of (n+1)-centered when binormal} yields
\begin{align*}
&\mbox{$\widetilde{T}_{\alpha,\beta}^{(n-1)}$ is binormal}\Longrightarrow \mbox{$\widetilde{T}_{\alpha,\beta}^{(n-1)}$ is $2$-centered},\\
&\mbox{$\widetilde{T}_{\alpha,\beta}^{(n-2)}$ is binormal, $\widetilde{T}_{\alpha,\beta}^{(n-1)}$ is $2$-centered}\Longrightarrow \mbox{ $\widetilde{T}_{\alpha,\beta}^{(n-2)}$ is $3$-centered},\\
&\mbox{$\widetilde{T}_{\alpha,\beta}^{(n-3)}$ is binormal, $\widetilde{T}_{\alpha,\beta}^{(n-2)}$ is $3$-centered}\Longrightarrow \mbox{ $\widetilde{T}_{\alpha,\beta}^{(n-3)}$ is $4$-centered},\\
&\centerline{\vdots}\\
&\mbox{$\widetilde{T}_{\alpha,\beta}^{(2)}$ is binormal, $\widetilde{T}_{\alpha,\beta}^{(3)}$ is $(n-2)$-centered}\Longrightarrow \mbox{ $\widetilde{T}_{\alpha,\beta}^{(2)}$ is $(n-1)$-centered},\\
&\mbox{$\widetilde{T}_{\alpha,\beta}$ is binormal, $\widetilde{T}_{\alpha,\beta}^{(2)}$ is $(n-1)$-centered}\Longrightarrow \mbox{ $\widetilde{T}_{\alpha,\beta}$ is $n$-centered},\\
&\mbox{$T$ is binormal, $\widetilde{T}_{\alpha,\beta}$ is $n$-centered}\Longrightarrow \mbox{ $T$ is $(n+1)$-centered}.\qedhere
\end{align*}
\end{proof}

\begin{remark}\label{rem:compare with Hilbert space operators-4} In the special case that $\alpha=\beta=\frac12$, the equivalence of items (i) and (ii) in the preceding corollary can be derived from \cite[Theorem~3.6]{Ito} and Lemma~\ref{lem:characterization of (n+1)-centered} for Hilbert space operators.
\end{remark}

 Suppose that $T\in\mathcal{L}(H)$ has the polar decomposition $T=U|T|$. Let
$\widetilde{U}$ be defined by \eqref{equ:defn of widetilde U}, and let  $\widetilde{U}^{(n)}$ be defined inductively  as

\begin{equation}\label{defn of U n}\widetilde{U}^{(0)}=U, \quad \widetilde{U}^{(1)}=\widetilde{U},\quad \widetilde{U}^{(n)}=\big[\widetilde{U}^{(n-1)}\big]^*\big[\widetilde{U}^{(n-1)}\big]^2\quad\mbox{for $n\ge 2$}.
\end{equation}
Specifically, if $T$ is $(n+1)$-centered, then since the operators $U^k (1\le k\le n+1$) are isometries, a simple calculation shows that
\begin{equation}\label{equ:computation result of U-(k)}
\widetilde{U}^{(k)}=(U^k)^*U^{k+1}\quad\mbox{for $1\leq k\leq n+1$}.
\end{equation}

\begin{theorem}\label{thm:characterization of n-centered by polar decomposition-T-alpha-beta} Suppose that $T\in\mathcal{L}(H)$ has the polar decomposition $T=U|T|$. Then for each $n\in\mathbb{N}$, the following statements are equivalent:
\begin{enumerate}
\item[{\rm (i)}] $T$ is $(n+1)$-centered;
\item[{\rm (ii)}] For every $\alpha>0$ and $\beta>0$, $\widetilde{T}_{\alpha,\beta}^{(k)}
    =\widetilde{U}^{(k)}\left|\widetilde{T}_{\alpha,\beta}^{(k)}\right|$ is the polar decomposition for $1\le k\le n$, in which $\widetilde{U}^{(k)}$ is given inductively by \eqref{defn of U n};
\item[{\rm (iii)}] There exist $\alpha>0$ and $\beta>0$ such that $\widetilde{T}_{\alpha,\beta}^{(k)}
    =\widetilde{U}^{(k)}\left|\widetilde{T}_{\alpha,\beta}^{(k)}\right|$ is the polar decomposition for $1\le k\le n$, in which $\widetilde{U}^{(k)}$ is given inductively by \eqref{defn of U n};.
\end{enumerate}
If items (i)--(iii) are satisfied, then  $\widetilde{U}^{(k)}$ is given by \eqref{equ:computation result of U-(k)} for  $1\leq k\leq n$.
\end{theorem}

\begin{proof} (i)$\Longrightarrow$(ii).  Assume that $T$ is $(n+1)$-centered. For every $\alpha>0$ and $\beta>0$, by Theorem~\ref{thm:binormality of T-alpha-beta-n}  $\widetilde{T}_{\alpha,\beta}^{(k)}$ is binormal for $0\le k\le n-1$. Let $\widetilde{U}^{(k)}$ be defined inductively by \eqref{defn of U n} for $0\le k\le n$. As is shown before,  $\widetilde{U}^{(k)}$ is given by \eqref{equ:computation result of U-(k)}.
Since $T=\widetilde{T}_{\alpha,\beta}^{(0)}$, which is binormal, so by
\cite[Theorem~4.14~(iii)]{Liu-Luo-Xu-1}
\begin{equation}\label{equ:polar decomposition of T-alpha-beta-1}
\widetilde{T}_{\alpha,\beta}^{(1)}
=\widetilde{U}^{(1)}\left|\widetilde{T}_{\alpha,\beta}^{(1)}\right|\quad  \mbox{is the polar decomposition}.
\end{equation}
Replacing $\widetilde{T}_{\alpha,\beta}^{(0)}$  with $\widetilde{T}_{\alpha,\beta}^{(1)}$, and using \cite[Theorem~4.14~(iii)]{Liu-Luo-Xu-1} once again, we see that
\begin{equation*}\label{equ:polar decomposition of T-alpha-beta-2}
\widetilde{T}_{\alpha,\beta}^{(2)}
=\widetilde{U}^{(2)}\left|\widetilde{T}_{\alpha,\beta}^{(2)}\right|\quad  \mbox{is the polar decomposition}.
\end{equation*}
Pursue  this process, we obtain eventually that
\begin{align*}\widetilde{T}_{\alpha,\beta}^{(n)}
    =\widetilde{U}^{(n)}\left|\widetilde{T}_{\alpha,\beta}^{(n)}\right|\quad\mbox{is the polar decomposition}.\end{align*}

The implication (ii)$\Longrightarrow$(iii) is clear.

(iii)$\Longrightarrow$(i). In view of \eqref{equ:polar decomposition of T-alpha-beta-1} and  (v)$\Longrightarrow$(ii) of \cite[Theorem~4.14]{Liu-Luo-Xu-1}, we conclude that $\widetilde{T}_{\alpha,\beta}^{(0)}$ is binormal. The same method employed shows that
 $\widetilde{T}_{\alpha,\beta}^{(1)}$, $\widetilde{T}_{\alpha,\beta}^{(2)}$, \dots, and $\widetilde{T}_{\alpha,\beta}^{(n-1)}$ are all binormal.
 Hence, by Theorem~\ref{thm:binormality of T-alpha-beta-n} $T$ is $(n+1)$-centered.
\end{proof}

We end this paper by a characterization of centered operators, which is immediate from
Theorems~\ref{thm:further describe of (n+2)-centered}, \ref{thm:cited lemma of Ito++}, \ref{thm:characterization of (n+1)-centered when binormal}, \ref{thm:binormality of T-alpha-beta-n},
\ref{thm:characterization of n-centered by polar decomposition-T-alpha-beta} and \cite[Corollary~4.5]{Liu-Luo-Xu-1}.

\begin{corollary}\label{cor:one result of centered operator-111} Let $T\in\mathcal{L}(H)$ have  the polar decomposition $T=U|T|$.  Then the following statements are equivalent:
\begin{enumerate}
\item[{\rm (i)}] $T$ is  centered;
\item[{\rm (ii)}] For every $\alpha>0$ and $\beta>0$, $\big[U^k|T|^{\alpha}(U^k)^*, |T|^{\beta}\big]=0$ for every $k\in\mathbb{N}$;
\item[{\rm (iii)}] There exist $\alpha>0$ and $\beta>0$ such that
$\big[U^k|T|^{\alpha}(U^k)^*, |T|^{\beta}\big]=0$ for every $k\in\mathbb{N}$;
\item[{\rm (iv)}] For every $\alpha>0$ and $\beta>0$, $\big[U^k|\widetilde{T}_{\alpha,\beta}|(U^k)^*, |\widetilde{T}_{\alpha,\beta}|\big]=0$ for every $k\in\mathbb{N}$;
\item[{\rm (v)}] There exist $\alpha>0$ and $\beta>0$ such that $\big[U^k|\widetilde{T}_{\alpha,\beta}|(U^k)^*, |\widetilde{T}_{\alpha,\beta}|\big]=0$ for every $k\in\mathbb{N}$;
\item[{\rm (vi)}] $T$ is binormal and for every $\alpha>0, \beta>0$, $\widetilde{T}_{\alpha,\beta}$ is centered;
\item[{\rm (vii)}] $T$ is binormal and there exist $\alpha>0$ and $\beta>0$ such that $\widetilde{T}_{\alpha,\beta}$ is centered;
\item[{\rm (viii)}] For every $\alpha>0$ and $\beta>0$, $\widetilde{T}_{\alpha,\beta}^{(k)}$ is $2$-centered for every $k\in \mathbb{N}\cup \{0\}$;
\item[{\rm (ix)}] There exist $\alpha>0$ and $\beta>0$ such that $\widetilde{T}_{\alpha,\beta}^{(k)}$ is $2$-centered for every $k\in \mathbb{N}\cup \{0\}$;
\item[{\rm (x)}] For every $\alpha>0$ and $\beta>0$, $\widetilde{T}_{\alpha,\beta}^{(k)}
    =\widetilde{U}^{(k)}\left|\widetilde{T}_{\alpha,\beta}^{(k)}\right|$ is the polar decomposition for every $k\in\mathbb{N}$, in which $\widetilde{U}^{(k)}$ is given inductively by \eqref{defn of U n};
\item[{\rm (xi)}] There exist $\alpha>0$ and $\beta>0$ such that $\widetilde{T}_{\alpha,\beta}^{(k)}
    =\widetilde{U}^{(k)}\left|\widetilde{T}_{\alpha,\beta}^{(k)}\right|$ is the polar decomposition for every $k\in\mathbb{N}$, in which $\widetilde{U}^{(k)}$ is given inductively by \eqref{defn of U n}.
\end{enumerate}
\end{corollary}

\begin{remark}\label{rem:compare with Hilbert space operators-5} In the special case that $\alpha=\beta=\frac12$,
 the equivalence of (i), (ix) and (xi) in the preceding corollary can be found in \cite[Theorem~4.1]{Ito} and \cite[Theorem~3.2]{Ito-Yamazaki-Yanagida} for Hilbert space operators.
 \end{remark}


\begin{thebibliography}{99}

\bibitem{Aluthge}
A. Aluthge, On $p$-hyponormal operators for $0<p<1$, \textit{Integr. Equ. Oper. Theory}, \textbf{13} (1990), 307--315.


\bibitem{Campbell-1}
S. L. Campbell, Linear operators for which $T^*T$ and $TT^*$ commute, \textit{Proc. Amer. Math. Soc.}, \textbf{34} (1972), 177--180.


\bibitem{Furuta}
T. Furuta, Generalized Aluthge transform on $p$-hyponormal operators, \textit{Proc. Amer. Math. Soc.}, \textbf{124} (1996), 3071--3075.


\bibitem{Ito}
M. Ito, T. Yamazaki, M. Yanagida, On the polar decomposition of the Aluthge transform and related results, \textit{J. Operator Theory}, \textbf{51} (2004), 303--319.


\bibitem{Ito-Yamazaki-Yanagida}
M. Ito, T. Yamazaki, M. Yanagida, On the polar decomposition of the product of two operators and its applications, \textit{Integr. Equ. Oper. Theory}, \textbf{49} (2004), 461--472.


\bibitem{Lance}
E. C. Lance, \textit{Hilbert $C^*$-modules--A toolkit for operator algebraists}, Cambridge University Press, Cambridge, 1995.


\bibitem{Liu-Luo-Xu}
N. Liu, W. Luo, Q. Xu, The polar decomposition for adjointable operators on Hilbert $C^*$-modules and centered operators, \textit{Adv. Oper. Theory}, \textbf{3} (2018), 855--867.


\bibitem{Liu-Luo-Xu-1}
N. Liu, W. Luo, Q. Xu, The polar decomposition for adjointable operators on Hilbert $C^*$-modules and $n$-centered operators, \textit{Banach J. Math. Anal.}, \textbf{13} (2019), 627--646.


\bibitem{Morrel}
B. B. Morre, P. S. Muhly, Centered operators, \textit{Studia Math.}, \textbf{51} (1974), 251--263.

\bibitem{Moslehian}M. S. Moslehian, Approximate $n$-idempotents and generalized Aluthge transform, \textit{Aequationes Math.}, \textbf{94} (2020), no. 5, 979--987.


\bibitem{XF}Q. Xu, X. Fang, A note on majorization and range inclusion of adjointable operators
on Hilbert $C^*$-modules, \textit{Linear Algebra Appl.}, \textbf{516} (2017), 118--125.


\end{thebibliography}
\end{document}